\documentclass[11pt]{amsart}
\usepackage{amsmath,amssymb,amsthm,upref,graphicx,mathrsfs, enumerate}

\usepackage{color,xcolor}

\numberwithin{equation}{section}
\allowdisplaybreaks

\newtheorem{theorem}{Theorem}[section]

\theoremstyle{definition}

\def\XXint#1#2#3{{\setbox0=\hbox{$#1{#2#3}{\int}$}
\vcenter{\hbox{$#2#3$}}\kern-.5\wd0}}

\usepackage{xcolor}
\numberwithin{equation}{section}
\allowdisplaybreaks[4]
\usepackage{amsmath,amsthm,amscd,amssymb,amsfonts, amsbsy}
\usepackage{latexsym}
\usepackage{exscale}

\def\XXint#1#2#3{{\setbox0=\hbox{$#1{#2#3}{\int}$}
\vcenter{\hbox{$#2#3$}}\kern-.5\wd0}}

\def\ch{\rm{ch}}
\def\dw{\textup{d}w}
\begin{document}

\title[$A_p$-$A_\infty$ estimates for general multilinear sparse operators]{$A_p$-$A_\infty$ estimates for general multilinear sparse operators}

\author{Mahdi Hormozi}
\address{University of Gothenburg \& Chalmers University of Technology SE-412 96 G\"oteborg, Sweden}
\email{hormozi@chalmers.se}

\author{Kangwei Li}
\address{Department of Mathematics and Statistics, P.O.B. 68 (Gustaf
H\"allstr\"omin katu 2b), FI-00014 University of Helsinki, Finland}
\email{kangwei.li@helsinki.fi}

\thanks{K.L. is supported by the European Union through T. Hyt\"onen's ERC
Starting Grant
`Analytic-probabilistic methods for borderline singular integrals'.
He is a member of the Finnish Centre of Excellence in Analysis and
Dynamics Research.
}
\date{\today}


\begin{abstract}
In this paper, we study the $A_p$-$A_\infty$ estimates for a class of multilinear dyadic positive operators. As applications, the $A_p$-$A_\infty$ estimates for different operators e.g. multilinear square functions and multilinear Fourier multipliers can be deduced very easily.
\end{abstract}
\maketitle

\section{Introduction}
The weighted norm inequality is a hot topic in harmonic analysis. In 1980s, Buckley \cite{Buckley} studied the quantitative relation between the weighted bound of Hardy-Littlewood maximal function and
the $A_p$ constant. Specifically, he showed that
 \[
 \|M\|_{L^p(w)}\le c [w]_{A_p}^{\frac 1 {p-1}},
 \]
 where recall that
 \[
 [w]_{A_p}:=\sup_{Q} \langle w\rangle_Q \langle w^{1-p'}\rangle_Q^{p-1}.
 \]
Here and through out, $\langle \cdot\rangle_Q$ denotes the average over $Q$.

 Since then,  the sharp weighted estimates for Calder\'on-Zygmund operators has attracted many authors' interest, which was referred to as the famous $A_2$ conjecture. The $A_2$ conjecture (now theorem) asserts that
 \[
 \|T\|_{L^p(w)}\le c [w]_{A_p}^{\max \{1,\frac 1{p-1}\}}.
 \]
It was finally proved by Hyt\"onen \cite{Hyt1}. The interested readers can consult \cite{Hyt2} for a survey on the history of the different proofs given for $A_2$ theorem. Moreover, Hyt\"onen and Lacey \cite{HL} extends the $A_2$ theorem to the so-called $A_p$-$A_\infty$ type estimates, i.e.,
\[
\| T(\cdot\sigma) \|_{L^p(\sigma)\rightarrow L^p(w)}\le c[w,\sigma]_{A_p}^{\frac 1p}([w]_{A_\infty}^{\frac 1{p'}}+[\sigma]_{A_\infty}^{\frac 1p}),
\]
where
\[
[w,\sigma]_{A_p}:= \sup_Q \langle w\rangle_Q \langle \sigma\rangle_Q^{p-1},\,\, [w]_{A_\infty}:=\sup_Q \frac 1{w(Q)}\int_Q  M(\mathbf 1_Q w)  dx
\]
and
 $\sigma$ needn't to be the dual weight of $w$, i.e., we don't require that $\sigma=w^{1-p'}$.

Now the story goes to the multilinear case. First we need to extend the $A_p$ weights to the multilinear case.
  Let $1<p_1,\dotsc,p_m<\infty$ and $p$ be numbers such that $\frac{1}{p}=\frac{1}{p_1}+\dotsb+\frac{1}{p_m}$ and denote $\vec  P = (p_1,\dotsc, p_m)$. Now we define  $[w,\vec \sigma]_{A_{\vec P}}$ constant:
\begin{equation*}
[w,\vec \sigma]_{A_{\vec P}}= \sup_Q\langle w\rangle_Q \prod_{i=1}^m \langle\sigma_i\rangle_Q^{\frac p{p_i'}}.
\end{equation*}
In the one weight case, i.e., $\sigma_i=w_i^{1-p_i'}$ and $w=\prod_{i=1}^m w_i^{p/{p_i}}$,  we say that $\vec{w}$ satisfies the $A_{\vec{P}}$ condition if $[w,\vec \sigma]_{A_{\vec P}}<\infty$, see \cite{LOPTT}.
  For the Buckley type estimate, the second author, Moen and Sun \cite{LMS} studied the sharp weighted estimates for multilinear maximal operators for all indices  and multilinear Calder\'on-Zygmund operators when $p>1$. The corresponding $A_p$-$A_\infty$ estimate was obtained in \cite{DLP} and \cite{LS-am}, respectively. Specifically, the result for multilinear maximal operators
  reads as
  \[
  \|\mathcal M(\cdot \vec\sigma)\|_{L^{p_1}(\sigma_1)\times\cdots \times L^{p_m}(\sigma_m)\rightarrow L^p(w)}\le [w,\vec \sigma]_{A_{\vec P}}^{\frac 1p}\prod_{i=1}^m [\sigma_i]_{A_\infty}^{\frac 1{p_i}}.
  \]
As to the multilinear Calder\'on-Zygmund operators, if $p>1$, then
\[
\|T(\cdot \vec\sigma)\|_{L^{p_1}(\sigma_1)\times\cdots \times L^{p_m}(\sigma_m)\rightarrow L^p(w)}\le [w,\vec \sigma]_{A_{\vec P}}^{\frac 1p}\Big(\prod_{i=1}^m [\sigma_i]_{A_\infty}^{\frac 1{p_i}}+[w]_{A_\infty}^{\frac 1{p'}} (\sum_{j=1}^m \prod_{i\neq j} [\sigma_i]_{A_\infty}^{\frac 1{p_i}}  )\Big).
\]

The spirit of the above results is reducing the problem to consider the so-called sparse operators. Recall that given a dyadic grid $\mathcal D$,  we say a collection $\mathcal S\subset \mathcal D$ is sparse if
\[
\Big|\bigcup_{\substack{Q'\subsetneq Q\\Q', Q\in \mathcal S}} Q'\Big|\le \frac 12 |Q|,
\]
and we denote $E_Q:= Q \setminus \cup_{Q'\in \mathcal S, Q'\subsetneq Q} Q'$.
Now given a sparse family $\mathcal{S}$ over a dyadic grid $D$ and $\gamma\geq 1$, a \textit{general multilinear sparse operator} is an averaging operator over $\mathcal{S}$ of the form

  \begin{equation*}
    T_{p_0,\gamma,\mathcal{S}}(\vec{f})(x)=\left(\sum_{Q\in S} \left[\prod_{i=1}^m \langle f_i \rangle_{Q,p_0}\right]^\gamma \chi_{Q}(x)\right)^{1/\gamma}
  \end{equation*}
where $p_0\in [1,\infty)$ and for any cube $Q$,
$$
\langle f \rangle_{Q,p_0} := \left(\frac{1}{|Q|} \int_Q |f(x)|^{p_0}dx \right)^{\frac{1}{p_0}}.
$$
It was proved in \cite{CR} that the multilinear Calder\'on-Zygmund operators are dominated pointwisely by $T_{1, 1,\mathcal S}$.  In \cite{BuiHormozi}, Bui and the first author also showed that the multilinear square functions are dominated pointwisely by $T_{1, 2,\mathcal S}$, and therefore, they obtained the Buckley type estimate for multilinear square functions.
For $\gamma=1$ and general $p_0$, it was shown in \cite{BCDH} that  $T_{p_0, 1,\mathcal S}$ can dominates a large class of operators with rough kernels (which include multilinear Fourier multipliers) as well. Therefore, everything are reduced to study $T_{p_0,\gamma,\mathcal{S}}$. Our main result states as follows.

        \begin{theorem}
        \label{main}
        Let $\gamma>0$.
 Suppose that $p_0<p_1,\dots,p_m<\infty $ with $\frac{1}{p}=\frac{1}{p_1}+\dots+\frac{1}{p_m}$. Let $w$ and $\vec\sigma$ be weights satisfying that $[w,\vec\sigma]_{A_{\vec P/p_0}}<\infty$ and $w, \sigma_i \in A_\infty$ for $i=1,\dots,m$. If  $\gamma \ge p_0$, then
\begin{align*}
\left\| T_{p_0,\gamma,\mathcal{S}}(\vec{f} )  \right\|_{L^p(w)}&\lesssim [w,\vec \sigma]_{A_{\vec P/p_0}}^{\frac 1p}\Big( \prod_{i=1}^m [\sigma_i]_{A_\infty}^{\frac 1{p_i}}+[w]_{A_\infty}^{(\frac 1\gamma-\frac 1p)_+}\sum_{j=1}^{m} \prod_{i\neq j}[\sigma_i]_{A_\infty}^{\frac 1{p_i}} \Big) \times \prod_{i=1}^m \|f_i\|_{L^{p_i}(w_i)},
\end{align*}
where $w_i=\sigma_i^{1-\frac{p_i}{p_0}}$, $i=1,\cdots,m$ and
\[
\left(\frac 1\gamma-\frac 1p\right)_+:=\max\left\{\frac 1\gamma-\frac 1p, 0\right\}.
\] If $\gamma < p_0$, then the above result still holds for all $p>\gamma $.
  \end{theorem}

The proof of Theorem \ref{main} is quite technical. In the literature, the $A_p$-$A_\infty$ estimates usually follows from testing condition. Our technique provide a way to obtain $A_p$-$A_\infty$ estimates without testing conditions. The idea follows from a recent paper by Lacey and the second author \cite {Lacey-Li}, where
  they studied the $A_p$-$A_\infty$ estimates for square functions in the linear case. We generalize their method to suit for the multilinear case with general parameters $\gamma$ and $p_0$.

\section{Proof of Theorem~\ref{main}}
Let us first observe that it suffices to prove Theorem~\ref{main} for $p_0=1$. Indeed, suppose Theorem~\ref{main} holds for $p_0=1$.  Consider the two weight norm inequality
\begin{equation}\label{eq:ep}
\|T_{p_0,\gamma,\mathcal S }(f  ,g  )\|_{L^p(w)}\le \mathcal N \|f\|_{L^{p_1}(w_1)} \|g\|_{L^{p_2}(w_2)},
\end{equation}
where we use $\mathcal N$ to denote the best constant such that \eqref{eq:ep} holds.
Rewrite \eqref{eq:ep} as
\[
\|T_{p_0,\gamma, \mathcal S}( f ^{1/{p_0}} , g ^{1/{p_0}} )\|_{L^p(w)}^{p_0}\le \mathcal N^{p_0} \|f^{1/{p_0}} \|_{L^{p_1}(w_1)} ^{p_0} \|g^{1/{p_0}} \|_{L^{p_2}(w_2)}^{p_0},
\]
which is equivalent to the following
\[
\|T_{1,\frac\gamma{p_0}, \mathcal S}( f   , g  )\|_{L^{p/{p_0}}(w)} \le \mathcal N^{p_0} \|f  \|_{L^{p_1/{p_0}}(w_1)}  \|g  \|_{L^{p_2/{p_0}}(w_2)}.
\]
Then by our assumption, we have
\[
\mathcal N \lesssim [w,\vec \sigma]_{A_{\vec P/{p_0}}}^{\frac 1p}\Big([\sigma_1]_{A_\infty}^{\frac 1{p_1}}[\sigma_2]_{A_\infty}^{\frac 1{p_2}}
 +[w]_{A_\infty}^{(\frac 1\gamma-\frac 1p)_+}([\sigma_1]_{A_\infty}^{\frac 1{p_1}}+[\sigma_2]_{A_\infty}^{\frac 1{p_2}})\Big).
\]

So we concentrate on the case $p_0=1$. As in \cite{LS-am}, we begin with $m=2$, that is we deal with the dyadic bilinear operators:
\[
T(f, g):=\Big( \sum_{Q\in\mathcal S} \langle f\rangle_Q ^\gamma\langle g\rangle_Q^\gamma \mathbf 1_Q  \Big)^{\frac 1\gamma}
\]
and we shall give the corresponding $A_p$-$A_\infty$ estimate.

Without loss of generality, we can assume that all cubes in $\mathcal S$ are contained in some root cube. As usual we only work on a subfamily $\mathcal S_a$, which is defined by the following
\[
\mathcal S_a:=\{Q \in \mathcal S: 2^a< \langle w\rangle_Q \langle\sigma_1\rangle_Q^{\frac p{p_1'}}\langle\sigma_2\rangle_Q^{\frac p{p_2'}}\le 2^{a+1}\}.
\]

Now we can define the principal cubes $\mathcal F$ for $(f,\sigma_1)$ and $\mathcal G$ for $(g, \sigma_2)$. Namely,
\begin{eqnarray*}
\mathcal F&:=& \bigcup_{k=0}^\infty \mathcal F_k, \quad \mathcal F_0:=  \{\textup{maximal cubes in }\mathcal S_a\}\\
\mathcal F_{k+1}&:=& \bigcup_{F\in \mathcal F_k}\ch_{\mathcal F}(F),\quad \ch_{\mathcal F}(F):= \{ Q\subsetneq F\, \textup{maximal \,s.t.} \langle f\rangle_Q^{\sigma_1}>2\langle f\rangle_F^{\sigma_1} \},
\end{eqnarray*}
and analogously for $\mathcal G$. We use $\pi_{\mathcal F}(Q)$ to denote the minimal cube in $\mathcal F$ which contains $Q$ and $\pi(Q)=(F,G)$ means that $\pi_{\mathcal F}(Q)=F$ and $\pi_{\mathcal G}(Q)=G$.
By our construction, it is easy to see that
\begin{equation}\label{eq:stopping}
\sum_{F\in \mathcal F} (\langle f\rangle_F^{\sigma_1})^{p_1}\sigma_1(F)\lesssim \|f\|_{L^{p_1}(\sigma_1)}^{p_1}.
\end{equation}

 We are going to prove that if $w,\sigma_1,\sigma_2$ be weights satisfying that $[w,\vec\sigma]_{A_{\vec P}}<\infty$ and $w,\sigma_1,\sigma_2\in A_\infty$.  Then
\begin{align*}
\|T(f\sigma_1, g\sigma_2)\|_{L^p(w)}&\lesssim [w,\vec \sigma]_{A_{\vec P}}^{\frac 1p}\Big([\sigma_1]_{A_\infty}^{\frac 1{p_1}}[\sigma_2]_{A_\infty}^{\frac 1{p_2}}\\
&+[w]_{A_\infty}^{(\frac 1\gamma-\frac 1p)_+}([\sigma_1]_{A_\infty}^{\frac 1{p_1}}+[\sigma_2]_{A_\infty}^{\frac 1{p_2}})\Big)\|f\|_{L^{p_1}(\sigma_1)}\|g\|_{L^{p_2}(\sigma_2)}.
\end{align*}

First, we consider the case $p\le \gamma $ with $\gamma \ge 1$.
\subsection{The case $p\le \gamma$ with $\gamma \ge 1$}
In this case, we have
\begin{align*}
&\Big\|  \Big( \sum_{Q\in\mathcal S_a} \langle f\sigma_1\rangle_Q ^\gamma\langle g\sigma_2\rangle_Q^\gamma \mathbf 1_Q  \Big)^{\frac 1\gamma}    \Big\|_{L^p(w)}\\
&=\Big\|  \Big( \sum_{Q\in\mathcal S_a} (\langle f\rangle_Q^{\sigma_1}) ^\gamma(\langle g\rangle_Q^{\sigma_2})^\gamma\langle\sigma_1\rangle_Q^\gamma\langle\sigma_2\rangle_Q^\gamma \mathbf 1_Q  \Big)^{\frac 1\gamma}    \Big\|_{L^p(w)}\\
&\lesssim \Big\|  \Big(\sum_{F\in\mathcal F}(\langle f\rangle_F^{\sigma_1}) ^\gamma \sum_{G\in\mathcal G}(\langle g\rangle_G^{\sigma_2})^\gamma \sum_{\substack{Q\in\mathcal S_a\\ \pi(Q)=(F,G)}}  \langle\sigma_1\rangle_Q^\gamma\langle\sigma_2\rangle_Q^\gamma \mathbf 1_Q  \Big)^{\frac 1\gamma}    \Big\|_{L^p(w)}\\
&\le   \Big(\sum_{F\in\mathcal F}(\langle f\rangle_F^{\sigma_1}) ^p \sum_{G\in\mathcal G}(\langle g\rangle_G^{\sigma_2})^p \Big\|\Big(\sum_{\substack{Q\in\mathcal S_a\\ \pi(Q)=(F,G)}}  \langle\sigma_1\rangle_Q^\gamma\langle\sigma_2\rangle_Q^\gamma \mathbf 1_Q  \Big)^{\frac 1\gamma} \Big\|_{L^p(w)}^p   \Big)^{\frac 1p}\\
&\lesssim  \Big(\sum_{F\in\mathcal F}(\langle f\rangle_F^{\sigma_1}) ^p \sum_{\substack{G\in\mathcal G\\ G\subset F}}(\langle g\rangle_G^{\sigma_2})^p \Big\| \sum_{\substack{Q\in\mathcal S_a\\ \pi(Q)=(F,G)}}  \langle\sigma_1\rangle_Q \langle\sigma_2\rangle_Q  \mathbf 1_Q    \Big\|_{L^p(w)}^p   \Big)^{\frac 1p}\\
&+  \Big( \sum_{G\in\mathcal G}(\langle g\rangle_G^{\sigma_2})^p \sum_{\substack{F\in\mathcal F\\ F\subset G}}(\langle f\rangle_F^{\sigma_1}) ^p \Big\| \sum_{\substack{Q\in\mathcal S_a\\ \pi(Q)=(F,G)}}  \langle\sigma_1\rangle_Q \langle\sigma_2\rangle_Q  \mathbf 1_Q    \Big\|_{L^p(w)}^p   \Big)^{\frac 1p}\\
&:=I+II.
\end{align*}
By symmetry we only focus on estimating $I$. From \cite{LS-am}, we already know that
\begin{align}\label{eq:ls}
&\Big\| \sum_{\substack{Q\in\mathcal S_a\\ \pi(Q)=(F,G)}}  \langle\sigma_1\rangle_Q \langle\sigma_2\rangle_Q  \mathbf 1_Q    \Big\|_{L^p(w)}\\
&\lesssim 2^{\frac ap}\Big(\sum_{\substack{Q\in\mathcal S_a\\ \pi(Q)=(F,G)}}\sigma_1 (Q)  \Big)^{\frac 1{p_1}}
\Big(\sum_{\substack{Q\in\mathcal S_a\\ \pi(Q)=(F,G)}}\sigma_2 (Q)  \Big)^{\frac 1{p_2}}\nonumber.
\end{align}
We also recall a fact that, for $\sigma\in A_\infty$ and $\mathcal S$ a sparse family, we have
\[
\sum_{\substack{Q\in\mathcal S\\ Q\subset R }} \sigma(Q) \le 2 \sum_{\substack{Q\in\mathcal S\\ Q\subset R }} \langle\sigma \rangle_Q |E_Q|\le 2\int_R M(\mathbf 1_R \sigma)dx\le 2[\sigma]_{A_\infty} \sigma(R).
\]
Therefore,
\begin{align*}
I &\lesssim 2^{\frac ap}    \Big(\sum_{F\in\mathcal F}(\langle f\rangle_F^{\sigma_1}) ^p \sum_{\substack{G\in\mathcal G\\ G\subset F}}(\langle g\rangle_G^{\sigma_2})^p     \Big(\sum_{\substack{Q\in\mathcal S_a\\ \pi(Q)=(F,G)}}\sigma_1 (Q)  \Big)^{\frac p{p_1}}
\Big(\sum_{\substack{Q\in\mathcal S_a\\ \pi(Q)=(F,G)}}\sigma_2 (Q)  \Big)^{\frac p{p_2}}           \Big)^{\frac 1p}\\
&\lesssim 2^{\frac ap} [\sigma_2]_{A_\infty}^{\frac 1{p_2}}   \Big(\sum_{F\in\mathcal F}(\langle f\rangle_F^{\sigma_1}) ^p \sum_{\substack{G\in\mathcal G\\ G\subset F}}(\langle g\rangle_G^{\sigma_2})^p     \Big(\sum_{\substack{Q\in\mathcal S_a\\ \pi(Q)=(F,G)}}\sigma_1 (Q)  \Big)^{\frac p{p_1}}
\sigma_2(G)^{\frac p{p_2}}           \Big)^{\frac 1p}\\
&\le 2^{\frac ap} [\sigma_2]_{A_\infty}^{\frac 1{p_2}}   \Big(\sum_{F\in\mathcal F}(\langle f\rangle_F^{\sigma_1}) ^p \Big(\sum_{\substack{G\in\mathcal G\\ \pi_{\mathcal F}(G)= F}}(\langle g\rangle_G^{\sigma_2})^{p_2} \sigma_2(G)\Big)^{\frac p{p_2}} \Big(\sum_{\substack{G\in\mathcal G\\ \pi_{\mathcal F}(G)= F}}\sum_{\substack{Q\in\mathcal S_a\\ \pi(Q)=(F,G)}}\sigma_1 (Q)  \Big)^{\frac p{p_1}}
          \Big)^{\frac 1p}\\
&\le 2^{\frac ap} [\sigma_2]_{A_\infty}^{\frac 1{p_2}}   \Big(\sum_{F\in\mathcal F}(\langle f\rangle_F^{\sigma_1}) ^{p_1} \sum_{\substack{G\in\mathcal G\\ \pi_{\mathcal F}(G)= F}}\sum_{\substack{Q\in\mathcal S_a\\ \pi(Q)=(F,G)}}\sigma_1 (Q)  \Big)^{\frac 1{p_1}}\\
&\times
\Big(\sum_{F\in \mathcal F}\sum_{\substack{G\in\mathcal G\\ \pi_{\mathcal F}(G)= F}}(\langle g\rangle_G^{\sigma_2})^{p_2} \sigma_2(G)\Big)^{\frac 1{p_2}}\\
&  \lesssim   2^{\frac ap}[\sigma_1]_{A_\infty}^{\frac 1{p_1}}  [\sigma_2]_{A_\infty}^{\frac 1{p_2}}  \|f\|_{L^{p_1}(\sigma_1)}\|g\|_{L^{p_2}(\sigma_2)},
\end{align*}
where \eqref{eq:stopping} is used in the last step.

\subsection{The case $p>\gamma $ with $p_1=\max\{p_1,p_2, q'\}$}
Here $q=p/\gamma$. By duality, we have
\begin{align*}
&\Big\|\Big( \sum_{Q\in \mathcal S_a}(\langle f\rangle_Q^{\sigma_1})^\gamma (\langle g\rangle_Q^{\sigma_2})^\gamma \langle\sigma_1\rangle_Q^\gamma\langle\sigma_2\rangle_Q^\gamma\mathbf 1_Q        \Big)^{\frac 1\gamma} \Big\|_{L^p(w)}^\gamma\\
&= \sup_{\|h\|_{L^{q'}(w)}=1}\sum_{Q\in \mathcal S_a}(\langle f\rangle_Q^{\sigma_1})^\gamma (\langle g\rangle_Q^{\sigma_2})^\gamma \langle h\rangle_Q^w \langle\sigma_1\rangle_Q^\gamma \langle\sigma_2\rangle_Q^\gamma w(Q).
\end{align*}
Now we suppress the supremum, and denote by $\mathcal H$ the principal cubes associated to $(h, w)$. Similarly, $\pi(Q)=(F,G,H)$ means that $\pi_{\mathcal F}(Q)=F$, $\pi_{\mathcal G}(Q)=G$ and $\pi_{\mathcal H}(Q)=H$. We have
\begin{align*}
&\sum_{Q\in \mathcal S_a}(\langle f\rangle_Q^{\sigma_1})^\gamma (\langle g\rangle_Q^{\sigma_2})^\gamma \langle h\rangle_Q^w \langle\sigma_1\rangle_Q^\gamma \langle\sigma_2\rangle_Q^\gamma w(Q)\\
&= \sum_{F\in\mathcal F} \sum_{\substack{G\in \mathcal G\\ G\subset F}} \sum_{\substack{H\in \mathcal H\\ H\subset G}}\sum_{\substack{Q \in\mathcal S_a\\ \pi(Q)=(F,G, H)}}(\langle f\rangle_Q^{\sigma_1})^\gamma (\langle g\rangle_Q^{\sigma_2})^\gamma \langle h\rangle_Q^w \langle\sigma_1\rangle_Q^\gamma \langle\sigma_2\rangle_Q^\gamma w(Q)\\
&+\sum_{G\in\mathcal G} \sum_{\substack{F\in \mathcal F\\ F\subset G}} \sum_{\substack{H\in \mathcal H\\ H\subset F}}\sum_{\substack{Q \in\mathcal S_a\\ \pi(Q)=(F,G, H)}}(\langle f\rangle_Q^{\sigma_1})^\gamma (\langle g\rangle_Q^{\sigma_2})^\gamma \langle h\rangle_Q^w \langle\sigma_1\rangle_Q^\gamma \langle\sigma_2\rangle_Q^\gamma w(Q)\\
&+\sum_{F\in\mathcal F} \sum_{\substack{H\in \mathcal H\\ H\subset F}} \sum_{\substack{G\in \mathcal G\\ G\subset H}}\sum_{\substack{Q \in\mathcal S_a\\ \pi(Q)=(F,G, H)}}(\langle f\rangle_Q^{\sigma_1})^\gamma (\langle g\rangle_Q^{\sigma_2})^\gamma \langle h\rangle_Q^w \langle\sigma_1\rangle_Q^\gamma \langle\sigma_2\rangle_Q^\gamma w(Q)\\
&+\sum_{G\in\mathcal G} \sum_{\substack{H\in \mathcal H\\ H\subset G}} \sum_{\substack{F\in \mathcal F\\ F\subset H}}\sum_{\substack{Q \in\mathcal S_a\\ \pi(Q)=(F,G, H)}}(\langle f\rangle_Q^{\sigma_1})^\gamma (\langle g\rangle_Q^{\sigma_2})^\gamma \langle h\rangle_Q^w \langle\sigma_1\rangle_Q^\gamma \langle\sigma_2\rangle_Q^\gamma w(Q)\\
&+ \sum_{H\in\mathcal H} \sum_{\substack{F\in \mathcal F\\ F\subset H}} \sum_{\substack{G\in \mathcal G\\ G\subset F}}\sum_{\substack{Q \in\mathcal S_a\\ \pi(Q)=(F,G, H)}}(\langle f\rangle_Q^{\sigma_1})^\gamma (\langle g\rangle_Q^{\sigma_2})^\gamma \langle h\rangle_Q^w \langle\sigma_1\rangle_Q^\gamma \langle\sigma_2\rangle_Q^\gamma w(Q)\\
&+ \sum_{H\in\mathcal H} \sum_{\substack{G\in \mathcal G\\ G\subset H}} \sum_{\substack{F\in \mathcal F\\ F\subset G}}\sum_{\substack{Q \in\mathcal S_a\\ \pi(Q)=(F,G, H)}}(\langle f\rangle_Q^{\sigma_1})^\gamma (\langle g\rangle_Q^{\sigma_2})^\gamma \langle h\rangle_Q^w \langle\sigma_1\rangle_Q^\gamma \langle\sigma_2\rangle_Q^\gamma w(Q)\\
&:= I+I'+II+II'+III+III'.
\end{align*}
First we estimate $I$. We have
\begin{align*}
I& \lesssim  \sum_{F\in\mathcal F}(\langle f\rangle_F^{\sigma_1})^\gamma \sum_{\substack{G\in \mathcal G\\ G\subset F}}(\langle g\rangle_G^{\sigma_2})^\gamma \sum_{\substack{H\in \mathcal H\\ H\subset G}}\langle h\rangle_H^w\sum_{\substack{Q \in\mathcal S_a\\ \pi(Q)=(F,G, H)}}   \langle\sigma_1\rangle_Q^\gamma \langle\sigma_2\rangle_Q^\gamma w(Q)\\
&\le \sum_{F\in\mathcal F}(\langle f\rangle_F^{\sigma_1})^\gamma \sum_{\substack{G\in \mathcal G\\ G\subset F}}(\langle g\rangle_G^{\sigma_2})^\gamma \int\Big(\sum_{\substack{H\in \mathcal H\\ H\subset G}} \sum_{\substack{Q \in\mathcal S_a\\ \pi(Q)=(F,G, H)}}   \langle\sigma_1\rangle_Q^\gamma \langle\sigma_2\rangle_Q^\gamma\mathbf 1_Q\Big) \Big(\sup_{\substack{H'\in \mathcal H\\ \pi_{\mathcal G} (H')= G}}\langle h\rangle_{H'}^w \mathbf 1_{H'}\Big)\dw\\
&\le  \sum_{F\in\mathcal F}(\langle f\rangle_F^{\sigma_1})^\gamma \sum_{\substack{G\in \mathcal G\\ G\subset F}}(\langle g\rangle_G^{\sigma_2})^\gamma
\Big\|     \sum_{\substack{H\in \mathcal H\\ H\subset G}} \sum_{\substack{Q \in\mathcal S_a\\ \pi(Q)=(F,G, H)}}   \langle\sigma_1\rangle_Q^\gamma \langle\sigma_2\rangle_Q^\gamma\mathbf 1_Q           \Big\|_{L^{q}(w)} \| \sup_{\substack{H'\in \mathcal H\\ \pi(H')= (F,G)}}\langle h\rangle_{H'}^w \mathbf 1_{H'} \|_{L^{q'}(w)}\\
&\le 2^{\frac {\gamma a} p} \sum_{F\in\mathcal F}(\langle f\rangle_F^{\sigma_1})^\gamma \sum_{\substack{G\in \mathcal G\\ G\subset F}}(\langle g\rangle_G^{\sigma_2})^\gamma
\Big(  \sum_{\substack{H\in \mathcal H\\ \pi(H)=(F,G)}} \sum_{\substack{Q\in \mathcal S_a\\ \pi(Q)=(F,G,H)}} \sigma_1(Q)\Big)^{\frac \gamma{p_1}}\\
&\times \Big(  \sum_{\substack{H\in \mathcal H\\ \pi(H)=(F,G)}} \sum_{\substack{Q\in \mathcal S_a\\ \pi(Q)=(F,G,H)}} \sigma_2(Q)\Big)^{\frac \gamma{p_2}}
\Big(  \sum_{\substack{H\in\mathcal H\\ \pi(H)=(F,G)}}      (\langle h\rangle_{H}^w)^{q'} w(H)            \Big)^{\frac 1{q'}}.
\end{align*}
Since $\frac \gamma{p_1}+\frac \gamma{p_2}+\frac 1{q'}=1$, by using H\"older's inequality twice we have
\begin{align*}
I &\lesssim 2^{\frac {\gamma a} p} \Big(  \sum_{F\in \mathcal F} (\langle f\rangle_F^{\sigma_1})^{p_1}  \sum_{\substack{G\in \mathcal G\\ \pi_{\mathcal F}(G)=F}}  \sum_{\substack{H\in \mathcal H\\ \pi(H)=(F,G)}} \sum_{\substack{Q\in \mathcal S_a\\ \pi(Q)=(F,G,H)}} \sigma_1(Q)\Big)^{\frac \gamma{p_1}}\\
&\times \Big(  \sum_{F\in \mathcal F}    \sum_{\substack{G\in\mathcal G\\ \pi_{\mathcal F}(G)=F}}    (\langle g\rangle_G^{\sigma_2})^{p_2}  [\sigma_2]_{A_\infty} \sigma_2(G)  \Big)^{\frac \gamma{p_2}}\Big(  \sum_{F\in\mathcal F}  \sum_{\substack{G\in \mathcal G\\ \pi_{\mathcal F}(G)=F}}  \sum_{\substack{H\in \mathcal H\\ \pi(H)=(F,G)}}    (\langle h\rangle_{H}^w)^{q'} w(H)     \Big)^{\frac 1{q'}}\\
&\overset{\eqref{eq:stopping}}{\lesssim} 2^{\frac {\gamma a} p} [\sigma_1]_{A_\infty}^{\frac \gamma{p_1}} [\sigma_2]_{A_\infty}^{\frac \gamma{p_2}} \|f\|_{L^{p_1}(\sigma_1)}^\gamma \|g\|_{L^{p_2}(\sigma_2)}^\gamma \|h\|_{L^{q'}(w)}.
\end{align*}
It is obvious that $I'$ can be estimated similarly. Next we estimate $II$.
We have
\begin{align*}
&\sum_{\substack{Q \in\mathcal S_a\\ \pi(Q)=(F,G, H)}}   \langle\sigma_1\rangle_Q^\gamma \langle\sigma_2\rangle_Q^\gamma w(Q)\\
&=\sum_{\substack{Q \in\mathcal S_a\\ \pi(Q)=(F,G, H)}}   \langle\sigma_1\rangle_Q^\gamma \langle\sigma_2\rangle_Q^\gamma \langle w\rangle_Q |Q|\\
&=\sum_{\substack{Q \in\mathcal S_a\\ \pi(Q)=(F,G, H)}}  (\langle w\rangle_Q)^{\frac {\gamma p_1'}{p}}(\langle\sigma_1\rangle_Q)^{\frac p{p_1'}\cdot \frac {\gamma p_1'}p}(\langle \sigma_2\rangle_Q)^{\frac p{p_2'}\cdot \frac{\gamma p_1'}p} \langle\sigma_2\rangle_Q^{\gamma-\frac{\gamma p_1'}{p_2'}}\langle w\rangle_Q^{1-\frac {\gamma p_1'}{p}}|Q|\\
&\lesssim 2^{\frac{\gamma p_1' a}p} \sum_{\substack{Q \in\mathcal S_a\\ \pi(Q)=(F,G, H)}}\langle\sigma_2\rangle_Q^{\gamma-\gamma\frac{p_1'}{p_2'}}\langle w\rangle_Q^{1-\frac{\gamma p_1'}p} |Q|
\end{align*}
Since $p_1=\max\{p_1,p_2,q'\}$ and $p>\gamma $, it is easy to check that
\[
0\le \gamma-\frac{\gamma p_1'}{p_2'}<1,\,\, 0\le 1-\frac{\gamma p_1'}p<1,
\]
and
\[
\frac 1r:=\gamma -\frac{\gamma p_1'}{p_2'}+ 1-\frac{\gamma p_1'}p<1.
\]
Set
\[
\frac 1s:= \gamma -\frac{\gamma p_1'}{p_2'} +\frac {1-\frac 1r}2.
\]
Then\[
\frac 1{s'}=1-\frac{\gamma p_1'}p  +\frac {1-\frac 1r}2,
\]
and therefore,
\begin{align*}
\sum_{\substack{Q \in\mathcal S_a\\ \pi(Q)=(F,G, H)}}   \langle\sigma_1\rangle_Q^\gamma \langle\sigma_2\rangle_Q^\gamma w(Q)
&\lesssim 2^{\frac{\gamma p_1' a}p} \int_G M(\sigma_2\mathbf 1_G)^{\gamma -\gamma\frac{p_1'}{p_2'}} M(w\mathbf 1_G)^{1-\frac{\gamma p_1'}p}dx\\
&\le  2^{\frac{\gamma p_1' a}p} \Big(\int_G M(\sigma_2\mathbf 1_G)^{s(\gamma-\gamma \frac{p_1'}{p_2'})} dx\Big)^{\frac 1s} \Big( \int_G M(w\mathbf 1_G)^{s'(1-\frac{\gamma p_1'}p)}\Big)^{\frac 1{s'}}
\end{align*}
Before we give further estimate, we introduce the {\em{Kolmogorov's inequality}} (see for example \cite{LOPTT}): Let $0<p<q< \infty$, then there exists a constant $C=C_{p,q}$ such that for any locally integrable function $f$,
$$
\| f \|_{L^p(Q, \frac{dx}{|Q|})} \leq C\| f \|_{L^{q,\infty}(Q, \frac{dx}{|Q|})}.
$$
With this inequality in hand, we have
\[
 \frac 1{|G|}\int_G M(w\mathbf 1_G)^{s'(1-\frac{\gamma p_1'}p)} dx
\le \| M(w\mathbf 1_G) \|_{L^{1,\infty}(G, \frac {dx}{|G|})}^{s'(1-\frac{\gamma p_1'}p)}\le \langle w\rangle_G^{s'(1-\frac{\gamma p_1'}p)},
\]
and
\[
\Big(\frac 1{|G|}\int_G M(\sigma_2\mathbf 1_G)^{s(\gamma -\gamma\frac{p_1'}{p_2'})} dx\Big)
\le \| M(\sigma_2\mathbf 1_G) \|_{L^{1,\infty}(G, \frac {dx}{|G|})}^{s(\gamma-\gamma\frac{p_1'}{p_2'})}\le \langle \sigma_2\rangle_G^{s(\gamma-\gamma\frac{p_1'}{p_2'})}.
\]
Thus we get
\begin{align*}
\sum_{\substack{Q \in\mathcal S_a\\ \pi(Q)=(F,G, H)}}   \langle\sigma_1\rangle_Q^\gamma \langle\sigma_2\rangle_Q^\gamma w(Q)
&\lesssim 2^{\frac{\gamma p_1' a}p}\langle \sigma_2\rangle_G^{\gamma -\gamma\frac{p_1'}{p_2'}}\langle w\rangle_G^{1-\frac{\gamma p_1'}p} |G|\\
&\lesssim 2^{\frac{\gamma a}p} w(G)^{1-\frac \gamma p} \sigma_1(G)^{\frac \gamma {p_1}}   \sigma_2(G)^{\frac \gamma {p_2}}.
\end{align*}
It follows that
\begin{align*}
II&\lesssim \sum_{F\in\mathcal F}(\langle f\rangle_F^{\sigma_1})^\gamma \sum_{\substack{H\in \mathcal H\\ H\subset F}} \langle h\rangle_H^w\sum_{\substack{G\in \mathcal G\\ G\subset H}}(\langle g\rangle_G^{\sigma_2})^\gamma\sum_{\substack{Q \in\mathcal S_a\\ \pi(Q)=(F,G, H)}}   \langle\sigma_1\rangle_Q^\gamma \langle\sigma_2\rangle_Q^\gamma w(Q)\\
&\lesssim \sum_{F\in\mathcal F}(\langle f\rangle_F^{\sigma_1})^\gamma \sum_{\substack{H\in \mathcal H\\ \pi_{\mathcal F}(H)= F}} \langle h\rangle_H^w\sum_{\substack{G\in \mathcal G\\ \pi_{\mathcal H}(G)= H}}(\langle g\rangle_G^{\sigma_2})^\gamma 2^{\frac{\gamma a}p} w(G)^{1-\frac \gamma p} \sigma_1(G)^{\frac \gamma {p_1}}   \sigma_2(G)^{\frac \gamma {p_2}}\\
&\lesssim 2^{\frac{\gamma a}p}[w]_{A_\infty}^{1-\frac \gamma {p}}[\sigma_1]_{A_\infty}^{\frac \gamma {p_1}}  \|f\|_{L^{p_1}(\sigma_1)}^\gamma \|g\|_{L^{p_2}(\sigma_2)}^\gamma \|h\|_{L^{q'}(w)},
\end{align*}
where again, the H\"older's inequality and \eqref{eq:stopping} are used in the last step.

Now we estimate $II'$. By similar arguments as that in the above, we have
\[
\sum_{\substack{Q \in\mathcal S_a\\ \pi(Q)=(F,G, H)}}   \langle\sigma_1\rangle_Q^\gamma \langle\sigma_2\rangle_Q^\gamma w(Q)
\lesssim 2^{\frac{\gamma a}p} w(F)^{1-\frac \gamma p} \sigma_1(F)^{\frac \gamma {p_1}}   \sigma_2(F)^{\frac \gamma {p_2}}.
\]
Then it follows that
\begin{align*}
II'&\lesssim \sum_{G\in\mathcal G}(\langle g\rangle_G^{\sigma_2})^\gamma \sum_{\substack{H\in \mathcal H\\ H\subset G}} \langle h\rangle_H^w\sum_{\substack{F\in \mathcal F\\ F\subset H}}(\langle f\rangle_F^{\sigma_1})^\gamma \sum_{\substack{Q \in\mathcal S_a\\ \pi(Q)=(F,G, H)}}   \langle\sigma_1\rangle_Q^\gamma \langle\sigma_2\rangle_Q^\gamma w(Q)\\
&\le
2^{\frac{\gamma a}p}[w]_{A_\infty}^{1-\frac \gamma{p}}[\sigma_2]_{A_\infty}^{\frac \gamma{p_2}}  \|f\|_{L^{p_1}(\sigma_1)}^\gamma \|g\|_{L^{p_2}(\sigma_2)}^\gamma \|h\|_{L^{q'}(w)}.
\end{align*}
$III$ and $III'$ can also be estimated similarly.

\subsection{The case $p>\gamma$ with $p_2=\max\{p_1,p_2, q'\}$}
By symmetry, this case can be estimated similarly as that in the previous subsection.

\subsection{The case $p>\gamma$ with $q'=\max\{p_1,p_2,q'\}$}
Again, we can decompose the summation to $I+I'+II+II'+III+III'$. The estimates of $I$ and $I'$ have no differences with the previous case. Now we consider $II$.
We have
\begin{align*}
\sum_{\substack{Q \in\mathcal S_a\\ \pi(Q)=(F,G, H)}}   \langle\sigma_1\rangle_Q^\gamma \langle\sigma_2\rangle_Q^\gamma w(Q)
&= \sum_{\substack{Q \in\mathcal S_a\\ \pi(Q)=(F,G, H)}}   \langle\sigma_1\rangle_Q^\gamma \langle\sigma_2\rangle_Q^\gamma \langle w\rangle_Q |Q|\\
&= 2^a \sum_{\substack{Q \in\mathcal S_a\\ \pi(Q)=(F,G, H)}}   \langle\sigma_1\rangle_Q^{\gamma-\frac p{p_1'}} \langle\sigma_2\rangle_Q^{\gamma-\frac p{p_2'}}   |Q|
\end{align*}
Since $p>\gamma$ and $q'\ge \max\{p_1,p_2\}$, we have
\[
\gamma- \frac p{p_1'}\ge0,\,\, \gamma-\frac p{p_2'}\ge 0,
\]
and
\[
\gamma- \frac p{p_1'}+ \gamma-\frac p{p_2'}=2\gamma+1-2p<1.
\]
Then follow the same arguments as that in the above, we get
\begin{align*}
\sum_{\substack{Q \in\mathcal S_a\\ \pi(Q)=(F,G, H)}}   \langle\sigma_1\rangle_Q^\gamma \langle\sigma_2\rangle_Q^\gamma w(Q)
&\lesssim 2^a  \langle\sigma_1\rangle_G^{\gamma-\frac p{p_1'}} \langle\sigma_2\rangle_G^{\gamma-\frac p{p_2'}}  |G|\\
&\lesssim 2^{\frac {\gamma a}p}\Big(\langle w\rangle_G\langle\sigma_1\rangle_G^{\frac p{p_1'}}\langle\sigma_2\rangle_G^{\frac p{p_2'}}\Big)^{1-\frac \gamma p}\langle\sigma_1\rangle_G^{\gamma-\frac p{p_1'}} \langle\sigma_2\rangle_G^{\gamma-\frac p{p_2'}}  |G|\\
&=2^{\frac{\gamma a}p} w(G)^{1-\frac \gamma p} \sigma_1(G)^{\frac \gamma {p_1}}   \sigma_2(G)^{\frac \gamma{p_2}}.
\end{align*}
Then follow the same arguments as the previous subsection we can get the desired conclusion. The estimates of $II'$, $III$ and $III'$ can also be estimated similarly.

\section{Applications}
 Theorem \ref{main} has some new applications. It is obvious that if an operator reduced to  $T_{p_0,\gamma,\mathcal{S}}$ for some $p_0$ and $\gamma$, then it is enough to apply  Theorem \ref{main} for those particular $p_0$ and $\gamma$. Thus, to find out the $A_p$-$A_\infty$ estimates for \textit{Multilinear square functions} (which were introduced and investigated in \cite{CXY, SXY, XY}), considering Proposition 4.2. of \cite{BuiHormozi}, it is enough to apply Theorem \ref{main} for $T_{1,2,\mathcal{S}}$ .\\

   To observe the other application, we first recall the class of multilinear {integral} operator which is bounded on certain products of Lebesgue spaces on $\mathbb R^n$ where associated kernel satisfies some mild regularity condition  which is weaker than the usual H\"older continuity of those in the class of multilinear Calder\'on-Zygmund singular integral operators. This class of the operators motivated from the recent works \cite{BD, GT1, KW,  LOPTT,  LMPR, LMRT,LRT} and weighted bounds for such operators studied in \cite{BCDH} very recently. The main example of such operators is \textit{Multilinear Fourier multipliers}. Now, to deduce the $A_p$-$A_\infty$ estimates for the operators of such class, it is enough to apply Theorem \ref{main} for $T_{p_0,1,\mathcal{S}}$ applying the main theorems of \cite{BCDH}. It is worth-mentioning that the $A_p$-$A_\infty$ estimates for linear Fourier multipliers was unknown as well as other noted multilinear operators.

\end{document}